# Quanto è importante risolvere e far risolvere problemi?

How important is to solve problems and to give problems to be solved?

Gilles Aldon

S2HEP, Université Lyon 1 Claude Bernard – Lione, Francia

gilles.aldon@ens-lyon.fr

> «Pour un esprit scientifique, toute connaissance est une réponse à une question. S'il n'y a pas eu de question, il ne peut y avoir connaissance scientifique. Rien ne va de soi. Rien n'est donné. Tout est construit»
>
> (Bachelard, 1967)

**Sunto** / Fare matematica implica tre livelli di manipolazione: manipolare l'astratto, manipolare i simboli e manipolare la logica. L'insegnamento della matematica passa quindi attraverso la proposta da parte del docente di situazioni in cui gli studenti possano esplorare una piccola parte della matematica attraverso queste manipolazioni. In questo modo, gli alunni lavorano sia sull'euristica che permette loro di confrontarsi con una vera e propria ricerca matematica sia sulle conoscenze in costruzione. Attraverso esempi di situazioni didattiche di problem solving, questo articolo intende mostrare come i problemi possano essere motori dell'apprendimento della matematica.

**Parole chiave:** problem solving; esperimenti; epistemologia; situazione didattica; apprendimento della matematica.

**Abstract** / Doing mathematics implies three levels of manipulation: manipulating the abstract, manipulating symbols and manipulating logic. Teaching mathematics therefore involves the teacher proposing situations in which pupils can explore a small part of mathematics through these manipulations. In doing so, the pupils work on both the heuristics enabling them to confront themselves with a real mathematical research and knowledge in construction. Through examples of didactic situations of problem solving, this article aims to show how problems can be drivers of mathematics learning.

**Résumé /** Faire des mathématiques, c'est manipuler l'abstrait, manipuler des symboles et manipuler la logique. Enseigner des mathématiques passe donc par la proposition par le professeur de situations dans lesquelles les élèves pourront à travers ces manipulations explorer une petite part des mathématiques. Ce faisant, les élèves travaillent à la fois les heuristiques permettant de se confronter à une véritable recherche mathématique et les connaissances en construction. A travers des exemples de situations didactiques de recherche de problèmes, cet article montre comment les problèmes peuvent être des moteurs de l'apprentissage des mathématiques.

**Mots-clefs:** résolution de problèmes; expériences; épistémologie; situation didactique; apprentissage des mathématiques.

## 1 Introduzione

Che cosa significa fare matematica? Questa domanda, all'apparenza semplice, apre un universo di quesiti di ordine filosofico, epistemologico e didattico. Stando alle definizioni date dal *Trésor de la Langue Française informatisé*,[1] il verbo "fare" utilizzato in questa domanda, ma anche utilizzato nel contesto più ampio del mondo dell'educazione matematica, suppone che il soggetto "darà origine alla", oppure "sarà l'autore della" matematica con la quale è confrontato. In questa prospettiva, "fare matematica" significa crearla. Tale risposta apparentemente semplice nasconde tuttavia temi centrali della filosofia della matematica, che hanno alimentato dibattiti fin dalla notte dei tempi: quando facciamo matematica stiamo scoprendo un mondo già esistente, come sostiene la filosofia di Platone, oppure stiamo costruendo la matematica affidandoci alla nostra intuizione di spazio e di tempo, come suggerisce Kant? Può anche darsi che queste riflessioni esulino dal quadro dell'insegnamento della matematica e abbiano senso unicamente per la ricerca, la creazione (o la scoperta!) matematica. Personalmente non lo credo. E questo articolo di riflessione tenterà di dimostrare che gli studenti, quando sono confrontati a dei problemi di matematica, diventano dei creatori di matematica, poiché manipolano dei *concetti a priori*: «i pensieri senza contenuto sono vuoti, le intuizioni senza concetti sono cieche» (Kant, 1781/1905, p. 91, *Introduzione alla Logica trascendentale*, traduzione dell'autore). Anche se da molto tempo le direttive istituzionali hanno messo l'accento sull'importanza dei problemi nell'insegnamento e nell'apprendimento della matematica,[2] l'imprecisione sui fondamenti epistemologici di tali direttive diventa un freno per la diffusione su grande scala di un insegnamento fondato sulla risoluzione dei problemi. In effetti, numerosi manuali o siti pedagogici si riferiscono a queste indicazioni per presentare una visione "scolastica" della risoluzione di problemi, appoggiandosi in particolare sulle ipotesi comportamentiste della costruzione delle conoscenze. Di conseguenza, per esempio, la risoluzione di un problema matematico a priori interessante e fecondo potrebbe essere ridotto al rispondere a una successione di domande collegate fra loro, il cui senso globale non potrà essere compreso se non attraverso uno sforzo supplementare di sintesi, sforzo richiesto molto raramente. A titolo d'esempio, si consideri un problema ampiamente studiato (Aldon et al., 2017); ecco riportati due possibili enunciati che si riferiscono alla stessa situazione matematica:

1. Siano *A* e *B* due punti situati nello stesso semipiano delimitato da una retta *d*. Dove occorre posizionare un punto *M* su *d*, in modo tale che la distanza *AM + MB* sia la minore possibile?
2. Siano *A* e *B* due punti situati nello stesso semipiano delimitato da una retta *d*. Sia *M* un punto di *d*.
    a. Costruire il punto *A'*, simmetrico al punto *A* rispetto a *d*.
    b. Dimostrare che *AM + MB = A'M + MB*.
    c. Dedurre la posizione di *M* su *d* in modo che la distanza *AM + MB* sia la minore possibile.

---

[1] http://atilf.atilf.fr/
[2] Mi riferisco in particolare al sistema educativo francese, i cui programmi scolastici sono disponibili in: https://www.ac-strasbourg.fr/fileadmin/pedagogie/mathematiques/Lycee/Programmes/programmes_Maths2019.pdf

Il primo enunciato lascia *aperta* la scelta della strategia da adottare: il metodo non è esplicitato e i mezzi di risoluzione sono quindi lasciati alla responsabilità degli allievi. Nel secondo enunciato, invece, la scelta del metodo, il quadro e gli strumenti da utilizzare sono imposti. Nonostante ciò, entrambi gli enunciati possono essere considerati come "problemi", a condizione che siano proposti a studenti per i quali rispondere all'enunciato non consista nello svolgere un esercizio d'applicazione, ma in una vera e propria ricerca e mobilitazione delle conoscenze pregresse da richiamare e riutilizzare. In particolare, il secondo enunciato è quello che in didattica della matematica (Arsac et al., 1991) si definisce un "problema aperto" per gli allievi, in analogia con i famosi sette "problemi del millennio", sei dei quali sono questioni ancora aperte che ancora non hanno trovato soluzione da parte dei matematici.

Le riflessioni proposte dal presente contributo evidenziano come la risoluzione di problemi permetta agli allievi di costruire e consolidare le loro conoscenze matematiche. Tali riflessioni, emerse durante l'implementazione in classe di quelle che chiameremo "situazioni didattiche di problem solving", poggiano sul lavoro di ricerca dell'equipe DREAM (*Démarche de Recherche pour l'Enseignement et l'Apprentissage des Mathématiques*) e sui risultati di una sperimentazione condotta nel quadro della *Maison des Mathématiques et de l'Informatique* (MMI) di Lione. Come verrà dettagliato in seguito, il lavoro dell'équipe DREAM si basa sull'insieme degli studi sviluppati attorno al tema del "problema aperto" in seno all'IREM di Lione[3] da più di trent'anni, nonché sulle ricerche condotte all'interno dell'Università di Lione (Arsac & al., 1991 ; Aldon & al., 2010). Una delle peculiarità del gruppo di lavoro è la sua composizione: ricercatori e docenti di ogni livello scolastico progettano e propongono delle sperimentazioni, riflettendo insieme sui loro risultati.

Nella prima parte di questo articolo, si delineano le ipotesi di lavoro dell'equipe DREAM attraverso alcune riflessioni preliminari che derivano dalla filosofia della matematica, che serviranno per precisare i fondamenti epistemologici alla base delle sperimentazioni presentate nella seconda parte. Tramite due esempi tratti dalle sperimentazioni svolte nel quadro di questa ricerca, mostrerò quali analisi sono possibili basando il lavoro in aula su queste ipotesi.

**2 Alcune riflessioni preliminari**

Giuseppe Longo (2020) ci ricorda che «la matematica è astratta, simbolica e rigorosa… al di là degli assiomi, dietro gli assiomi» (p. 3, traduzione dell'autore). È proprio riprendendo queste tre direzioni di pensiero (astratta, simbolica, logica) che presenterò la posizione dell'equipe DREAM sul tema dei problemi nell'insegnamento della matematica (Aldon et al., 2017; Front, 2015; Gardes, 2013).

*2.1 Che cos'è un problema?*

Da un punto di vista etimologico, la parola "problema" deriva dal greco *Πρόβλημα,* il cui verbo di riferimento è *προβάλλειν*, ossia "lanciare in avanti": si tratta quindi di "proiettarsi in avanti", di dirigersi verso il futuro, dunque di tentare di utilizzare le conoscenze per risolvere una questione, utilizzando dei metodi riconosciuti dalla società, dalla disciplina o più generalmente dall'istituzione nella quale evoluiamo. In particolare, in questo articolo, si assume la prospettiva di Brousseau secondo il quale:

> «Tutto ciò rende non solo vero, ma anche interessante un teorema e, insieme a questo, ciò che Gonseth chiamava il carattere idoneo di una conoscenza matematica, ossia ciò che fa sì che questa conoscenza esista come soluzione ottimale [...] [è] una soluzione a un problema»
> 
> (Brousseau, 1976, p. 103, traduzione dell'autore).

Così, all'interno della matematica, ossia nel processo di matematizzazione verticale (Treffers, 1987) su cui è focalizzato questo articolo, un problema è una domanda che conduce a un risultato matematico, un teorema, che amplierà il campo delle conoscenze della disciplina. Risolvere un problema implica quindi l'utilizzo di ragionamenti riconosciuti come validi all'interno della disciplina. Questo ci conduce evidentemente a una prima domanda: quali sono questi ragionamenti validi? Domanda, questa, che ci riporta all'apprendimento di questi ragionamenti e delle conoscenze necessarie per arrivare alla soluzione del problema con il quale siamo confrontati. Queste conoscenze sono in effetti una condizione necessaria per far sì che il problema possa essere affrontato. Schoenfeld (1998) mostra come sia necessario attivare, nella risoluzione di problemi, sia conoscenze dichiarative (1 e 2) sia conoscenze procedurali (3), che consistono in:

1. le conoscenze informali, che possono essere associate alle cosiddette *intuizioni pure* di Kant, perché si fondano sulle forme di sensibilità dello spazio e del tempo;
2. i fatti e le definizioni, che permetteranno di nominare e manipolare gli oggetti matematici in gioco nel problema;
3. le procedure, che si riferiscono sia alle regole di calcolo sia ai ragionamenti e alle routine che permettono di manipolare gli oggetti e di alimentare e rafforzare così la deduzione e la validazione dei risultati.

---

[3] Institut de Recherche sur l'Enseignement des Mathématiques de Lyon: http://math.univ-lyon1.fr/irem/

Sotto queste ipotesi, risolvere un problema di matematica consiste nel manipolare simultaneamente *l'astratto*, i *simboli* e la *logica*, come sottolineato da Bonnay e Dubucs (2011): «[possiamo] dare una rappresentazione fedele e cognitivamente plausibile dei tre elementi che sono al cuore della matematica: primo, gli oggetti a cui il matematico si riferisce, secondo, le formule che usa, e terzo, la sua attività mentale» (pp. 17-18, traduzione dell'autore).

## *2.2 Manipolare l'astratto*

Il titolo di questo paragrafo potrebbe sembrare contraddittorio, soprattutto se si considera che una "manipolazione" consiste nel tenere "in mano" l'oggetto che vogliamo trasformare o comprendere. È proprio qui che si pone tutta la questione dell'"esistenza" degli oggetti matematici. Se adottiamo un punto di vista platonico, il mondo matematico esiste indipendentemente dall'uomo; si tratta di scoprire le proprietà di oggetti che esistono a priori. Ma questo realismo si riferisce a oggetti matematici che esistono indipendentemente dall'uomo o alla verità degli enunciati matematici?

Stiamo parlando di quello che Petitot (1986) definisce «la differenza ontologica tra fenomeno e oggetto», per cui la definizione di oggetto diventa «il semplice operatore di traducibilità dei fatti empirici e/o dei dati numerici in un linguaggio formale» (p. 3, traduzione dell'autore). È quindi sulla base di questa differenza ontologica che il rapporto tra matematica e realtà oggettiva può essere descritto in termini dialettici: «è in una prospettiva trascendentale che è possibile interpretare al meglio questo investimento dell'astratto nella genesi del concreto» (Lautman, 1935-1939, citato da Petitot, 1986, p. 14, traduzione dell'autore); così la matematica si situa all'incrocio tra l'esperienza dei fenomeni e il mondo delle idee. È ciò che permette questa manipolazione degli oggetti a condurre alla creazione fruttuosa di concetti matematici. A questo punto occorre essere prudenti e non confondere l'esperienza primitiva, che consiste in una manipolazione non finalizzata alla costruzione di conoscenza, e l'esperienza scientifica, che ha senso solo quando vi è una riflessione sui risultati dell'esperienza. La verità scientifica emerge allora come il risultato di un doppio movimento: l'*inter*-azione con il concreto e la *trans*-azione nel costruire l'esperienza e nell'interpretare i suoi risultati. I due movimenti partecipano dialetticamente all'elaborazione dei fatti.

*2.3 Manipolare i simboli*

La sintassi matematica è essenziale per la costruzione del ragionamento, ma la sistematizzazione del formalismo matematico conduce inevitabilmente ad allontanarsi dal senso e, quindi, dal rapporto con le esperienze fenomenologiche. Le tesi di Frege e Russel, che nel XIX e XX secolo hanno tentato di ridurre la matematica alla logica, si sono scontrate con il senso e con le realtà dei fenomeni: «L'astratto non può rivendicare un'esistenza autonoma: questo è sufficiente per confutare l'idea secondo la quale gli assiomi rappresentano delle convenzioni poste liberamente dalla mente» (Gonseth, 1936/1974, p. 92, traduzione dell'autore).

D'altro canto, fare matematica implica l'utilizzo del simbolismo logico e dei sistemi di segni che permettono di costruire e comunicare la matematica. È proprio una problematica dell'insegnamento quella di potere padroneggiare i diversi sistemi di segni e di potere comunicare passando dall'uno all'altro (Duval, 1991, 1995). Gli oggetti matematici sono infatti manipolati attraverso le loro rappresentazioni all'interno di molteplici sistemi di segni, dove tutto il problema consiste nel comprendere che l'oggetto stesso sarà caratterizzato dall'insieme delle sue rappresentazioni; in questo modo, un oggetto matematico può essere definito come la classe di equivalenza di tutte le sue rappresentazioni. Questo fatto pone due conseguenze essenziali:

1. un oggetto matematico può essere intelligibile in un contesto e sconosciuto o difficile in un altro;
2. la conversione da un registro di rappresentazione ad un altro è essenziale per concepire un oggetto matematico, cosa che richiede un lavoro di traduzione che talvolta fa perdere degli elementi di significato e altre volte ne aggiunge di nuovi; modificando il significante, cioè la maniera di designare l'oggetto, si modifica, si arricchisce e si completa il significato, ossia l'oggetto designato.

Per illustrare il primo punto, consideriamo ad esempio questa moltiplicazione:

$$2375 \times 7 = 16625$$

Abbiamo numerosi modi per valutare la correttezza, o quantomeno la plausibilità, di questo risultato: l'ordine di grandezza (2500 moltiplicato per 7 dà come risultato 17500, l'ordine di grandezza è rispettato), la cifra delle unità ($7 \times 5 = 35$ e il risultato deve terminare con 5), e inoltre la moltiplicazione in sé non è difficile e la verifica può essere svolta facilmente. Consideriamo ora questa moltiplicazione:

$$2375 \times 7 = 21353$$

Il nostro primo riflesso potrebbe essere quello di dire che il risultato è falso, utilizzando gli stessi criteri di verifica mobilitati per la prima operazione. Tuttavia, se cambiamo contesto di partenza e consideriamo che la moltiplicazione è fatta in base 8, ecco che improvvisamente ci sentiamo smarriti, e la verifica della validità del risultato non è così scontata: le tabelline della moltiplicazione in base 8 non ci sono familiari e $7 \times 5 = 43$ non suona bene al nostro orecchio! Eppure, stiamo manipolando sempre dei numeri naturali, solo attraverso un'altra loro rappresentazione.

La seconda osservazione riguarda la semiotica e ha numerose conseguenze didattiche che sono state ampiamente studiate in letteratura (si vedano, ad esempio, Arzarello et al., 2009; Bartolini Bussi et al., 2005; Duval, 1991, 1995; Hitt, 2004; Presmeg et al. 2018; Sabena, 2018). A questo proposito, Duval sottolinea che: «Le rappresentazioni semiotiche sono delle produzioni costituite da segni appartenenti ad un sistema di rappresentazione che possiede i propri vincoli di significato e di funzionamento» (Duval, 1991, p. 234, traduzione dell'autore). In seguito, Duval (1993) insiste sul legame tra le conversioni di registri e l'apprendimento: «Non ci può essere un vero apprendimento finché le situazioni e i compiti proposti non prendono in conto la necessità di diversi registri di rappresentazione, per il funzionamento cognitivo del pensiero e la centralità dell'attività di conversione» (p. 64, traduzione dell'autore).

Il rapporto dialettico tra ~~noésis~~ *néosis* (comprensione concettuale di un oggetto) e *sémiosis* (comprensione delle rappresentazioni semiotiche di un oggetto) è al centro della comprensione degli oggetti della matematica e partecipa all' «implicazione dell'astratto nella genesi del concreto» di cui parla Lautman (1977, p. 205, traduzione dell'autore).

In quest'ottica, anche la manipolazione dei simboli, essenziale in ogni pratica matematica, contribuisce in questo modo alla comprensione dell'oggetto manipolato in un rapporto dialettico la cui importanza in termini didattici è spesso minimizzata, in particolare lasciando intendere che un oggetto sia intimamente legato a una delle sue rappresentazioni. L'esempio precedente illustra come nel campo numerico questa difficoltà sia presente e, anche se un obiettivo essenziale della scuola primaria è la padronanza del sistema di numerazione decimale e posizionale, è illusorio e talvolta pericoloso lasciare intendere che un numero sia equivalente alla sua rappresentazione in questo sistema di segni.

## 2.4 Manipolare la logica

Una delle prime osservazioni di Polya (1945/1957) è la seguente: «se non riesci a risolvere il problema proposto, prova a risolvere prima i problemi ad esso correlati» (p. XVII, traduzione dell'autore); un consiglio essenziale per invitare gli allievi ad osare, proponendo dei "piccoli frammenti di matematica" che potranno forse non condurre a una soluzione generale, ma che, ad ogni modo, parteciperanno alla costruzione delle conoscenze matematiche dell'allievo.

*2.4.1 Il ruolo della logica nell'interpretazione di un'esperienza*

Quando noi risolviamo un problema – dove il "noi" rappresenta ogni persona confrontata con un problema matematico (sia esso un matematico o una matematica, un allievo o un professore) – tutti i mezzi sono consentiti. È esattamente in questa fase che la creatività e l'aiuto dell'intuizione permettono di mettere in evidenza degli elementi essenziali per la risoluzione. In questo senso, la risoluzione di un problema è il luogo in cui poter realizzare delle esperienze matematiche. A proposito di esperienza, soffermiamoci un momento sul senso che diamo a questo termine in matematica: seguendo quanto espresso da Dias (2008), «l'andirivieni tra teoria ed esperienza è precisamente quello che caratterizza una pratica di tipo sperimentale» (p. 27, traduzione dell'autore); la sperimentazione ha senso unicamente quando si articolano la formulazione di ipotesi e la (o le) loro validazione(i) che può (o possono) intendersi come una verifica empirica o una dimostrazione nel senso matematico del termine (tornerò sulla nozione di validazione nel prossimo paragrafo). Nella sperimentazione si suppone che il soggetto possa manipolare degli oggetti concreti o sufficientemente familiari perché possano apparire come tali. Di conseguenza, un allievo di scuola media potrà manipolare i numeri naturali e le operazioni elementari su questi numeri *come se* si trattasse di oggetti reali. L'espressione "come se" rimanda a delle posizioni epistemologiche forti, che ritroviamo per esempio nel modo in cui Poincaré (1902/1968) utilizza il concetto di etere, ripreso in seguito da Mizony (2006):

> «[…] un etere è la reificazione (la cosificazione) di uno spazio matematico utilizzato per studiare un campo di fenomeni. E se il campo di fenomeni è unico, vi è una molteplicità di spazi matematici in grado di esprimere un campo della fisica (quello che Poincaré definisce come il pluralismo teorico), e quindi una molteplicità di eteri possibili se reifichiamo questi spazi matematici».
>
> (p. 92, traduzione dell'autore)

L'esperienza può essere realizzata anche senza utilizzare le logiche della matematica. Tuttavia, vi è un momento in cui se si intende trasmettere i concetti compresi o intuiti attraverso l'esperienza occorre passare attraverso le regole ammesse dalla logica matematica che sono al centro dell'enunciato di una frase matematica. È passando al vaglio della logica le frasi così prodotte che possono essere validati o rigettati i risultati dell'esperienza empirica. In effetti, la deduzione logica esiste eccome in matematica e fonda i ragionamenti "ipotetico-deduttivi", ma il senso della matematica dimora nella sua relazione con il reale. Per esempio, la retta euclidea senza spessore è una pura astrazione e fonda tutta la geometria che si applica intuitivamente al mondo che ci circonda.

Se dunque da un lato l'esperienza fonda la scoperta della matematica, dall'altro questa scoperta deve tenere conto della struttura logica intrinseca all'esperienza: «L'attività matematica è un'attività sperimentale, in altre parole è un sistema di atti legalizzati da regole e sottomesso a delle condizioni che non dipendono da esse» (Petitot, 1987, p. 98, traduzione dell'autore) citando la posizione di Cavaillès nel dibattito del 4 febbraio 1939 alla Società Francese di Filosofia (Cavaillès & Lautman, 1945). O, come direbbe Lautman, la matematica opera dialetticamente un passaggio dell'essenza all'esistenza:

> «Passiamo in maniera impercettibile dalla comprensione di un problema dialettico alla genesi di un universo di nozioni matematiche. È al riconoscimento di questo momento, in cui l'Idea dà origine al reale che deve, a mio avviso, mirare alla Filosofia matematica».
>
> (Lautman, 1977, p. 147, traduzione dell'autore)

### 2.4.2 La fase di validazione

Come abbiamo detto pocanzi, la validazione è un momento fondamentale della risoluzione di un problema, che permette all'andirivieni tra esperienza, teoria e riflessione sui risultati dell'esperienza. La risoluzione di un problema può essere concepita come una situazione didattica (Brousseau, 1986) che deve essere devoluta agli allievi. In una situazione, Brousseau (1986) classifica le interazioni dell'alunno con il *milieu* (ovvero, in senso ampio, l'ambiente di apprendimento) in tre grandi categorie, che corrispondono a tre fasi:
- *fase di azione*, nella quale avvengono scambi di informazioni non codificate in un linguaggio, che corrispondono a delle azioni che i protagonisti fanno direttamente sul *milieu* e sugli altri protagonisti della situazione, interpretandone le retroazioni e attivando le proprie conoscenze e i teoremi-in-atto di cui dispongono;
- *fase di formulazione*, nella quale avvengono scambi di informazioni codificate in un linguaggio;
- *fase di validazione*, nella quale avvengono scambi di giudizi.

Nella fase di validazione avviene la messa in relazione tra i risultati dell'esperienza sul *milieu* e le conoscenze dei protagonisti; è in questa fase che le regole della logica differiscono dalle regole instaurate dalla logica classica, che sottintendono il ragionamento ipotetico-deduttivo della matematica. Brousseau scrive:

> «Le dimostrazioni e le validazioni esplicite dovrebbero appoggiarsi le une sulle altre fino a giungere all'evidenza, ma la loro articolazione non è per forza automatica. I saperi e le conoscenze si attualizzano in un'attività di ricerca o di prova secondo le modalità che l'euristica cerca di scoprire e che l'intelligenza artificiale tenta di riprodurre».
>
> (Brousseau, 1986, p. 348, traduzione dell'autore)

In questa fase si tratta allora di costruire un messaggio matematico preciso, "vero" nel senso matematico del termine, in una dimensione dialogica il cui oggetto verte prevalentemente sulla veridicità delle asserzioni.

Ora, in questo andirivieni tra l'esperienza e la teoria matematica che dovrebbe modellizzarla, il senso degli oggetti manipolati è primordiale. In una fase di validazione si tratta quindi di decidere da una parte sulla verità semantica degli enunciati, dall'altra sulla verità pragmatica mettendo in gioco il linguaggio naturale. Per garantire la veridicità degli enunciati e permettere l'avanzamento della risoluzione di un problema, i protagonisti devono fare affidamento su una logica, diversa dalla logica matematica delle proposizioni, ma sufficientemente esplicita perché possa permettere di fare avanzare la trasformazione dei confronti avuti con il *milieu* verso una traduzione formale, matematicamente valida nel linguaggio formale della matematica, e alla quale è possibile attribuire il valore di verità "vero" nella logica delle proposizioni. Per essere più precisi, occorre distinguere tra il credere che un'asserzione sia vera ("*véracité*") e l'accordo tra i protagonisti della situazione circa il fatto che l'asserzione è vera ("*véridicité*", intesa come accordo sulla "*véracité*" di una asserzione) (Vernant, 2004, 2008).

La modellizzazione della logica dialogica è allora costruita su un'articolazione tra *véridicité* e verità, vale a dire una coerenza di dialogo che possa essere giudicata internamente dai protagonisti, e un risultato che possa essere confrontato con un giudizio esterno. La verità può allora intendersi, come proposto da Tarski, come relativa a un meta-linguaggio nel quale la verità è definita da:

*P* è vera se e solamente se *p*, dove *p* è la proposizione espressa da *P*.

Questa definizione semantica della verità concede delle interpretazioni diverse e da un punto di vista filosofico è neutra:

> «Infatti, la definizione semantica della verità non implica nulla riguardo alle condizioni per le quali un enunciato come "la neve è bianca" possa essere dichiarato. Implica solamente che ogni volta che asseriamo o rifiutiamo questo enunciato, dobbiamo essere pronti ad asserire o a rifiutare l'enunciato correlato: "L'enunciato 'la neve è bianca' è vero"» (Tarski, 1972, p. 295, traduzione dell'autore).

Così, la fase di validazione proposta da Brousseau può essere modellizzata in questa prospettiva, come sottolineato da Durand-Guerrier:

«Diciamo che l'apparato logico di cui abbiamo bisogno per trattare la questione dell'apprendimento della dimostrazione e del ragionamento nella prospettiva, in particolare, della teoria delle situazioni didattiche, deve essere più consapevole dell'attività effettiva del matematico, e pensiamo abbiano dato delle prove di fecondità, per questo, i concetti e i metodi della teoria dei modelli di Tarski».

(Durand-Guerrier, 2005, p. 23, traduzione dell'autore)

L'insieme di queste considerazioni filosofiche e logiche costituiscono un fondamento per analisi delle situazioni didattiche di problem solving proposte dall'equipe DREAM. L'originalità di queste situazioni consiste nel lasciare agli allievi la scelta di decidere le esperienze da realizzare, i tipi di ragionamento da attivare e le teorie matematiche alle quali attingere. Nella seconda parte dell'articolo utilizzeremo le considerazioni appena fatte in questo paragrafo per analizzare il lavoro di allievi confrontati a due situazioni di questo tipo.

**3 Due esempi di situazioni didattiche di problem solving**

Gli esempi che seguono sono tratti dalle osservazioni condotte dall'equipe DREAM e nel contesto di una sperimentazione condotta alla MMI in una scuola superiore e in alcune scuole elementari e medie della regione di Lione (Aldon & Garreau, 2017). Il primo esempio esposto verterà sull'analisi dei comportamenti di allievi di *Première*, ovvero, in Francia, del secondo anno di liceo (grado 11). Il secondo esempio illustra i comportamenti di allievi di *Cours Moyen 2* (ultimo anno di scuola elementare, grado 5) e di *Sixième* (primo anno di scuola media, grado 6).

*3.1 Primo esempio: le frazioni egizie*

*3.1.1 Contesto*

Oltre che sul problema aperto, i lavori dell'equipe DREAM si basano anche sugli studi sviluppati attorno all'articolazione tra logica e ragionamento matematico (Durand-Guerrier, 2005), e le tesi di Marie-Line Gardes (2013) e di Mathias Front (2015). La metodologia di ricerca si basa sul paradigma della *design-based research* per la quale la ricerca è fondata su una necessità epistemologica per i ricercatori di agire con gli insegnanti (Monod-Ansaldi et al., 2019; Nizet et al.,

2019). In questo approccio, partendo da situazioni matematiche ricche, l'equipe analizza dal punto di vista matematico, didattico e pragmatico il potenziale di queste situazioni al fine di trasformarle in situazioni didattiche di problem solving. Oltre che ai lavori accademici (master e tesi di dottorato), le pubblicazioni dell'equipe di ricerca (si vedano per esempio: Aldon et al., 2010; Aldon et al., 2012; Front, 2012; Front & Gardes, 2015) mettono in evidenza gli apprendimenti degli allievi confrontati a delle situazioni didattiche di problem solving. Ognuno dei problemi analizzati viene sperimentato in alcune classi, ed è proprio sull'analisi di tali osservazioni che si costruisce il prossimo paragrafo. Il problema presentato nel seguito è stato proposto a differenti livelli scolastici, ma nell'articolo mi riferirò unicamente all'osservazione svolta in una classe *Première* (grado 11).

*3.1.2 Metodologia*

L'obiettivo del lavoro dell'equipe DREAM è quello di proporre agli insegnanti di matematica degli approcci per utilizzare i problemi nel loro insegnamento. Le seguenti domande sono al centro della ricerca:

- Quali sono le conoscenze, le competenze trasversali e meta-matematiche che è possibile valutare in una pratica di problem solving? E quali sono gli indicatori che possono essere messi in atto?
- La creatività e l'inventiva matematica sviluppate nelle situazioni didattiche di problem solving modificano l'immagine che hanno gli studenti della matematica (e il loro desiderio di fare matematica)? E per gli insegnanti?
- Le situazioni didattiche di problem solving che sviluppano una forma di acquisizione della conoscenza aiutano gli studenti a progredire in altre aree dell'attività matematica? Come gli studenti reinvestono in altri contesti le competenze e le conoscenze sviluppate?

Per rispondere a queste domande, stiamo realizzando nelle classi degli insegnanti coinvolti nel progetto delle sperimentazioni che vengono osservate e analizzate alla luce delle ipotesi di ricerca basate sulle nostre posizioni epistemologiche descritte nella prima parte di questo articolo (**par. 2**). Le osservazioni si sono svolte in aula con un osservatore per gruppo che ha registrato i dialoghi degli allievi.

*3.1.3 Il problema*

Così, è stata organizzata una sperimentazione in una classe di liceo con l'obiettivo di raccogliere elementi di risposta alle domande precedenti, in particolare per quanto riguarda la creatività e l'invenzione matematica sviluppata da questo tipo di problemi negli studenti di una classe *Première* di liceo scientifico. L'analisi seguente si concentra sull'adeguatezza delle nostre ipotesi epistemologiche e didattiche confrontate con la realtà di una situazione didattica di problem solving in aula. In questo esempio, i dialoghi degli studenti che hanno lavorato in gruppi formati dall'insegnante sono stati registrati e l'analisi si concentra sia sul problema in sé – il suo potenziale per sviluppare la creatività e l'apprendimento, e il suo legame con le conoscenze matematiche degli studenti –, sia sul confronto tra le basi epistemologiche del nostro studio e il lavoro effettivo degli studenti.

Questo esempio è stato scelto per evidenziare due aspetti sviluppati nel **par. 2**:

- da un lato, la manipolazione dei simboli senza la comprensione del fenomeno in questione porta all'impossibilità di concludere;
- d'altra parte, il tentativo di una costruzione di una verità matematica si fa attraverso la ricerca di un accordo tra i partecipanti alle tappe del ragionamento.

L'enunciato che è stato proposto dall'insegnante in questa classe *Première* è il seguente:

1. Riesci a trovare due numeri interi naturali *a* e *b* distinti, in modo tale che: $1 = 1/a + 1/b$?
2. Riesci a trovare tre numeri interi naturali *a*, *b* e *c* distinti in modo tale che: $1 = 1/a + 1/b + 1/c$?
3. Riesci a trovare quattro numeri interi naturali *a*, *b*, *c* e *d* distinti in modo tale che: $1 = 1/a + 1/b + 1/c + 1/d$?

Continua…

Da un punto di vista matematico, alla prima domanda si può rispondere negativamente considerando diverse strategie, tra cui per esempio sfruttando la decrescita della funzione inversa su N:

Possiamo supporre che $2 \leq a < b$ e quindi $1/b < 1/a \leq 1/2$ di conseguenza $1/a + 1/b < 1$.

La seconda domanda può essere affrontata trovando un metodo generale per determinare una scomposizione in *n* frazioni dell'unità. Per esempio, partendo dalle uguaglianze $1 = 1/2 + 1/2$ e $1/2 = 1/3 + 1/6$, possiamo dedurre che $1 = 1/2 + 1/3 + 1/6$.

Possiamo individuare il generarsi di una soluzione di ordine superiore, per esempio osservando che qualsiasi sia $n > 0$, $1/n = 1/(n+1) + 1/n(n+1)$; da questo segue per esempio che:

$1 = 1/2 + 1/3 + 1/7 + 1/42$ scomponendo $1/6$ in $1/7 + 1/42$. E così via.

La lezione si è svolta in un'ora, gli allievi sono stati raggruppati in sette gruppi ciascuno dei quali è stato osservato e registrato da osservatori esterni. L'insegnante ha esplicitato dall'inizio che gli allievi avevano diritto di passare alla seconda domanda anche se non avevano risolto la prima.

L'esperienza dimostra in effetti che la prima domanda può essere difficile per degli allievi a causa della risposta negativa: non esistono due interi naturali *a* e *b* distinti tale che $1=1/a+1/b$! Questa risposta negativa da un lato può entrare in contrasto con il contratto abituale della classe (quando una domanda è posta, vi è sempre una soluzione), dall'altro si scontra con il fatto che la risposta negativa implica una dimostrazione *universale*: per ogni coppia di numeri interi *a* e *b* distinti, $1 \neq 1/a+1/b$, mentre una risposta positiva implicherebbe una dimostrazione *esistenziale*: esistono due numeri interi naturali *a* e *b* distinti tali che: $1=1/a+1/b$.

La consegna data dall'insegnante è quella di risolvere in gruppo il problema e di prendere nota su un cartellone, prima della fine dell'ora, dei risultati trovati da ciascun gruppo. La richiesta da parte dell'insegnante di scrivere un resoconto sul cartellone ha lo scopo di stimolare la validazione dei ragionamenti, di permettere in un primo tempo agli allievi di capirsi tra loro e, conseguentemente, di andare verso una verità semantica espressa nel cartellone.

Non mi soffermerò nel dettagliare l'insieme delle produzioni degli allievi, ma mi attarderò piuttosto su alcuni momenti e dialoghi significativi in termini di manipolazione e di validazione all'interno dei gruppi.

*3.1.4 Osservazione e discussione*

Esaminiamo in particolare questo dialogo tenutosi in un gruppo di quattro allievi (A1, B1, C1 e D1). Dopo un piccolo momento di riflessione individuale, il dialogo inizia così.

1. A1: Hai risolto la domanda 2?
2. D1: Sì.
3. A1: Una sola soluzione? E se moltiplichi tutto per due? Prova con i multipli.

[…]

4. A1: Bisogna trovare un numero intero che divida 1 e che sia inferiore a *b*. Bisognerebbe avere *a<b* oppure *b<a*.
5. B1: Non puoi scegliere.
6. A1: È una questione di multipli.
7. B1: Di multipli?

[…]

8. D1: La frazione più piccola con dei numeri interi è 1/2; non puoi avere una frazione di interi più piccola [Mima con le mani l'unità, poi la metà e mostra che resta una metà che non potrà essere colmata da una delle frazioni egizie successive].
9. C1: Ah sì, hai ragione.
10. D1: La domanda 1 non è possibile, la frazione più piccola è ½.

11. B1: Puoi spiegarti?

12. D1: Per la domanda 2, la soluzione è: *a = 2, b = 3, c = 6*; e per il 3 è: *a = 2, b = 3, c = 9, d = 18*.

13. A1: Perché la domanda 1 non è possibile?

14. D1 [si lancia allora nella determinazione di somme di 4, 5 frazioni egizie che costruisce scomponendo l'ultimo termine in due frazioni egizie e che spiega nel modo seguente ai suoi compagni]: Ne sono sicuro. Ho preso la cifra più piccola senza arrivare a 1. È il trucco delle frazioni egizie: 1-1/2 = 1/2, 1/2-1/3=1/6, ho preso il numero sopra: 1/6-1/9=1/18. Si può continuare mettendo al posto del 18…

15. A1: Perché la prima domanda è impossibile?

Questo estratto è particolarmente interessante per via di un'incomprensione nata tra i protagonisti. D1, all'inizio silenzioso, entra nella discussione avendo capito perfettamente l'impossibilità di scomporre di 1 nella somma di due frazioni egizie distinte. La giustificazione che propone si basa sulla manipolazione delle grandezze (intervento 8); il suo gesto mostra chiaramente la sua comprensione profonda del fenomeno, che l'allievo traduce in maniera più confusa oralmente: «la frazione più piccola con due numeri interi è 1/2»; i suoi compagni, che non si trovano nel suo stesso registro di grandezze ma si situano piuttosto nel registro algebrico, come si vede in seguito, ritengono che questa spiegazione orale non sia convincente, poiché con dei numeri interi maggiori di 2 tutte le frazioni egizie sono più piccole di 1/2. La confusione nel discorso di D1 tra «il più piccolo numero intero positivo che possiamo considerare è 2» e «la frazione più piccola è 1/2» (intervento 10), non gli permette di comunicare efficacemente il suo ragionamento ai compagni, come testimoniato dalle domande di B1 (intervento 11) e di A1 (intervento 15). Il «Ah sì, hai ragione!» di C1 (intervento 9) potrebbe farci supporre che D1 sia riuscito almeno a convincere C1. Tuttavia, quando si tratta di riassumere la loro ricerca sul cartellone, C1 scrive:

$$1) \quad \frac{1}{a} + \frac{1}{b} = 1 \quad donc \quad \frac{1b + 1a}{ab} = 1$$

$$\frac{a+b}{ab} = 1 \quad donc \quad a+b = ab \quad \text{impossible!}$$

**Figura 1**. Giustificazione per la prima domanda da parte di C1.

C1 ritorna quindi a un ragionamento algebrico in cui l'ultimo passaggio (cfr. «*Impossible*» in **Figura 1**) non è giustificato. Per C1, la manipolazione dei simboli nel campo algebrico è consolidata, ma il suo confronto con D1 (interventi 8-10) non si trasforma in una traduzione formale, perché risultano diversi i domini in cui lavorano D1 da una parte e C1 dall'altra (insieme anche a A1 e B1).

Il secondo aspetto interessante di questo estratto consiste invece nella manipolazione dei registri di rappresentazione che D1 effettua, traduce e trasforma in scrittura formale. Tutto si svolge come se il fatto di avere compreso e gestito la situazione di due frazioni egizie gli permettesse di cambiare registro e di tradurre nel registro aritmetico e in scrittura formale adeguata i risultati per trovare finalmente una maniera di generare delle serie di frazioni egizie, lunghe quanto desidera, di somma 1. È quello che D1 detta a C1 nella seconda parte del cartellone (anche se mentre ricopia C1 commette un errore scrivendo al punto 2 in **Figura 2**, 1-1/3=1/6 al posto di 1/2-1/3=1/6).

**Figura 2.** Generazione di una serie di frazioni egizie di somma 1.

Riassumendo, in questo gruppo vediamo bene che le manipolazioni dell'astratto, dei simboli e della logica dipendono fortemente dal registro utilizzato; il passaggio dall'uno all'altro non può avvenire senza un riferimento al campo fenomenologico che soggiace a questo registro.

Esaminiamo ora il comportamento di un secondo gruppo nel quale, dopo diversi tipi di manipolazione, le esperienze successive sfociano in una dimostrazione formale quasi completa. In un primo momento il gruppo di quattro allievi, che chiameremo A2, B2, C2 e D2, si situano in un registro algebrico e deducono che se $1=1/a+1/b$ allora $a=b-1/a$ e $b=a-1/b$. In seguito, rimpiazzando *a* con il suo valore nella seconda uguaglianza, trovano 1=1, che non sanno bene come interpretare. Di conseguenza cambiano strategia e sperimentano con i numeri.

1. C2: $1=1/2+1/2$ dunque è possibile, ma *a=b*. E poi c'è confusione.
2. A2: Dev'esserci un algoritmo per passare dalla prima domanda alla seconda, bisogna quindi risolvere la prima.
   [Malgrado una discussione sulle frazioni, non arrivano a convincersi dell'impossibilità della domanda 1]

3. A2: Una soluzione c'è, ma non so come spiegartela!

[Dopo essere passato alla seconda domanda, ritorna sulla prima]

4. A2: 1/2 + 1/3, poi 1/3+1/4, ci si allontana da 1.

5. C2: Ma non stiamo facendo le osservazioni di scienze della vita e della Terra!

6. A2: Sì, serve vedere.

[…]

7. A2: Io ho provato che più è grande più ci si allontana.

Infine, scrivono il loro ragionamento sul cartellone in **Figura 3**.

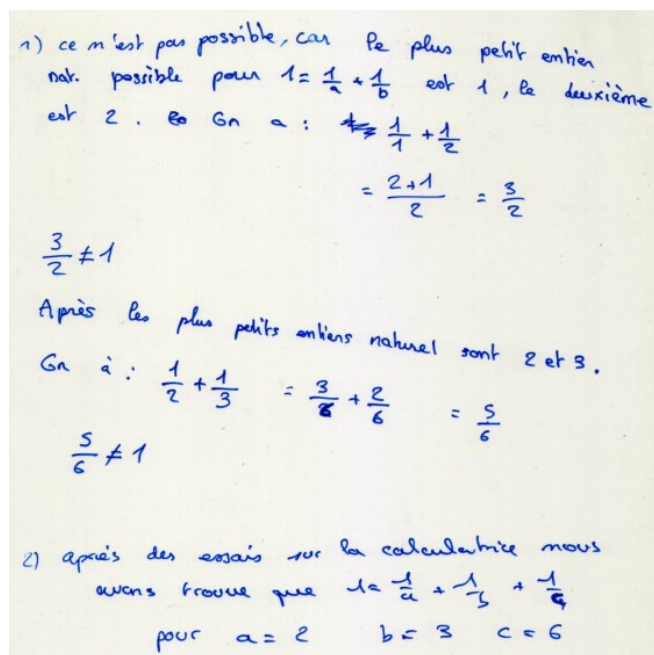

**Figura 3**. Il cartellone completo del secondo gruppo.

In questo caso, anche se la scrittura formale non è completamente raggiunta, il ragionamento seguito dal gruppo conduce effettivamente a una dimostrazione dell'impossibilità di scomporre 1 in una somma di due frazioni egizie distinte. È proprio attraverso un dialogo che verte sulla veridicità delle asserzioni che viene elaborata la costruzione matematica del ragionamento.

In questo gruppo il passaggio da una verità pragmatica, fondata sulla manipolazione di esempi, a una verità semantica espressa nel cartellone, si traduce nel confronto tra l'orale e lo scritto. Il risultato ottenuto dal gruppo va a smontare l'ipotesi iniziale formulata da A2: «Deve esserci un algoritmo per passare dalla prima domanda alla seconda, bisogna quindi risolvere la prima» (intervento 2). Da notare anche che, nonostante le sollecitazioni dell'insegnante, questa ipotesi resta presente, a dimostrazione di quanto sia difficile scardinare alcune clausole del contratto didattico. Dopo avere trovato una terna in risposta alla seconda domanda gli allievi ritornano sulla prima, ma questa volta lo fanno per confutare la congettura iniziale e dimostrarne l'impossibilità utilizzando in maniera pragmatica la decrescita della funzione $x \to 1/x$ («Io ho provato che più [il denominatore] è grande, più ci si allontana [da 1]», dice A2).

Inoltre, è interessante il breve dialogo tra A2 e C2 sull'esperienza in matematica a confronto con quella condotta nella disciplina scienze della vita e della Terra: «Sì, serve vedere» dice A2, situando in questo modo l'esperienza e la manipolazione dell'astratto matematico rispetto all'esperienza su oggetti concreti effettuata nelle scienze sperimentali.

In generale, le osservazioni di tutti i gruppi mettono in evidenza la difficoltà di uscire da un'investigazione algebrica, che è indotta dall'enunciato dato sotto questa forma. Un numero importante di allievi si lancia in tentativi di dimostrazione utilizzando il calcolo algebrico, in generale senza successo. Da ciò, potremmo ipotizzare una confusione tra l'utilizzo delle lettere e l'algebra, provocata senza dubbio dall'insegnamento dell'algebra sovente presentata come calcolo "con le lettere". Il problema che si pone qui è un problema di aritmetica e una risoluzione della prima domanda potrebbe essere espressa dal seguente ragionamento per assurdo:

> Supponiamo che esistano $a$ e $b$ tali che $1/a+1/b=1$
> ne deduciamo allora che $a+b=ab$
> quindi $a$ divide $a+b$ e quindi $a$ divide $b$ e allo stesso modo $b$ divide $a+b$ e quindi $b$ divide $a$.
> $a=kb=kk'a$ quindi $kk'=1$ e quindi $k=k'=1$, da cui $a=b$, cosa che contraddice l'ipotesi. Quindi si può concludere che non esistono due numeri naturali $a$ e $b$ distinti tali che $1/a+1/b=1$.

Ma questo ragionamento suppone che la congettura sia inizialmente prodotta, spostando peraltro il problema in un dominio poco lavorato in classe.

### *3.2 Secondo esempio: Il problema del foglio di carta*

*3.2.1 Contesto*

La *Maison des Mathématiques et de l'Informatique* (MMI) è un luogo «di mediazione dei saperi dedicato alle scienze matematiche e informatiche attraverso un approccio vivo, ludico e pluridisciplinare» (https://www.mmi-lyon.fr/, traduzione dell'autore). Si tratta anche di un centro di risorse pedagogiche il cui obiettivo è quello di accompagnare gli insegnanti e i loro allievi nella diffusione di metodi d'insegnamento e apprendimento, mettendo in evidenza l'attività degli allievi e la scoperta. Nel contesto delle attività della MMI, ho partecipato alla formazione di insegnanti di scuola elementare e di scuola media con l'obiettivo di sensibilizzarli all'utilizzo di problemi matematici nel loro insegnamento. In questo quadro, e in relazione con gli istituti, ho proposto una ricerca collaborativa che si è svolta per diversi anni consecutivi, raggruppando allievi di *CM2* (grado 5, ultimo anno di scuola elementare) e allievi di *Sixième* (grado 6, primo anno di scuola media). La finalità di questo dispositivo era duplice: da una parte si voleva coinvolgere gli insegnanti in un procedimento in cui i problemi di matematica fossero al cuore del loro insegnamento, dall'altra si voleva dare la possibilità agli allievi di vivere una situazione matematica non scolastica, invitandoli a condurre una riflessione sulle proprie procedure di risoluzione. Si è cercato di coinvolgere gli allievi in un lavoro che conducesse alla ricerca di modalità risolutive partendo da un problema sufficientemente complesso perché la ricerca potesse svilupparsi nel tempo, arricchendo le lezioni di matematica sia da un punto di vista delle procedure utilizzate che dei concetti matematici manipolati. Ho così proposto a tutte le classi di risolvere il seguente problema. Tutto parte da un foglio di carta che taglieremo in più pezzi.

> Immaginiamo di tagliarlo in due, poi di prendere uno dei due pezzi e di tagliarlo in due, poi di prendere uno dei pezzetti e di tagliarlo in due e così di seguito. Quante volte devo svolgere questa operazione per ottenere 2016 pezzi di carta?
> Ora, lo taglio in tre, poi prendo uno dei tre pezzi e lo taglio in tre, poi prendo uno dei pezzetti e lo taglio in tre e così di seguito. Potrei riuscire ad ottenere 2016 pezzi?[4]
> E se facessi la stessa operazione, però tagliando il foglio in quattro pezzi? E in cinque?... Più in generale, quanti sono i tagli che mi permetteranno di ottenere 2016 pezzi? E se volessi ottenerne 2017? E 2018?

---

[4] Questo problema è stato sperimentato nel 2016! Fatto interessante anche da un punto di vista matematico poiché i divisori di 2015 sono {1, 5, 13, 31, 65, 155, 403, 2015} e quindi i valori di $n$ che permettono di ottenere 2016 sono: {2, 6, 14, 66, 156, 404, 2016}!

*3.2.2 Metodologia*

Le osservazioni del lavoro in classe realizzate in questo caso differiscono dalle osservazioni realizzate nella prima sperimentazione descritta nel **par. 3.1**. In questo caso, infatti, le osservazioni vertono sulla continuità del lavoro e si appoggiano sulle produzioni degli allievi, talvolta scritte, altre volte orali, durante la presentazione dei loro lavori.

La sperimentazione si è svolta come segue: sono state coinvolte otto classi, quattro *CM2* e quattro *Sixième*. All'inizio, ho presentato il problema alle classi tutte raccolte in un anfiteatro e gli alunni hanno potuto realizzare la loro ricerca iniziale in questo quadro. In un secondo tempo, gli insegnanti si sono messi in gioco nel supervisionare le ricerche degli allievi durante i due mesi della sperimentazione, organizzando nelle loro classi dei momenti di ricerca delle soluzioni al problema chiedendo agli alunni di fare periodicamente una valutazione della loro ricerca. Nella terza fase, ho fatto visita alle otto classi in modo che gli studenti potessero presentarmi le loro ricerche. Questi momenti di lavoro sono stati l'occasione per me di ascoltare gli alunni, prendere nota delle loro relazioni e rilanciare la ricerca. Infine, il quarto momento ha visto un raggruppamento delle otto classi durante il quale gli alunni hanno presentato gli uni agli altri lo stato delle loro ricerche che avevano accuratamente riassunto su dei cartelloni. L'analisi di questo lavoro si basa quindi sulle varie tracce raccolte durante il processo di ricerca e sulle interviste condotte con gli insegnanti. Questo esempio è stato scelto per evidenziare la relazione tra la manipolazione di oggetti concreti (in questo caso, pezzi di carta) e la manipolazione di oggetti astratti (numeri) che porta alla manipolazione di regole logiche attraverso la costruzione e la manipolazione delle rappresentazioni simboliche dei concetti coinvolti. Così, le tre manipolazioni descritte nel **par. 2** contribuiscono alla comprensione del problema e alla sua risoluzione, almeno parzialmente.

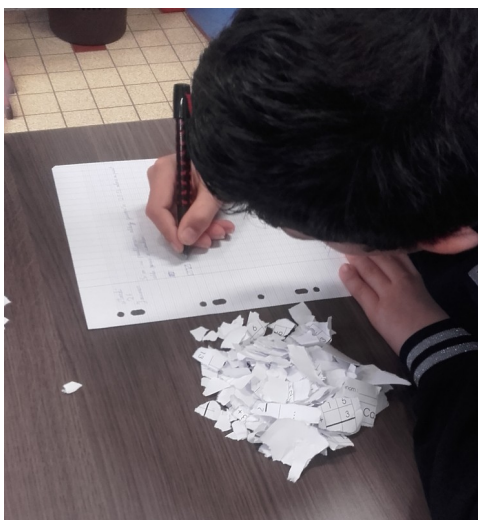

**Figura 4**. Il passaggio dall'esperienza alla riflessione sull'esperienza.

Chiaramente, le esperienze concrete hanno costituito la prima fase della ricerca di soluzione al problema e la prima domanda (tagliare in due parti) è stata risolta empiricamente: all'inizio c'è un foglio, quindi un pezzo di carta; poi due, poi tre e via di seguito. Tutti i numeri sono raggiunti e, in particolare, 2016 sarà raggiunto alla duemilasedicesima tappa.[5] Questa prima esperienza è semplice e contribuisce così alla devoluzione del problema. Numerosi insegnanti mi hanno riferito come questa tappa fosse stata particolarmente importante per molti allievi, in particolare per coloro che presentavano delle difficoltà matematiche; il problema si lascia affrontare! Ma le difficoltà iniziano quando il problema diventa più complesso, ed è in questo momento che l'andirivieni tra l'esperienza e i concetti soggiacenti struttura la comprensione del fenomeno, permettendo di rinforzare le conoscenze. Accade proprio questo quando il numero di ritagli del foglio aumenta! Fare l'esperienza concreta diventa presto impossibile e la scelta di chiedere se 2016 (o comunque un numero sufficientemente grande) sia raggiungibile obbliga l'allievo ad abbandonare l'esperienza concreta per passare a manipolare gli oggetti matematici che permettono di modellizzare l'esperienza (**Figura 4**). Di conseguenza, il passaggio alla riflessione, che è una condizione necessaria per operare un distacco da oggetti concreti verso il trattamento di oggetti matematici, è una tappa fondamentale per il processo sperimentale in matematica.

In tal modo, nel procedimento di risoluzione di questo problema, gli allievi hanno progressivamente costruito la formalizzazione del problema, come attestato da questo riassunto (**Figura 5**) di un allievo di *Sixième* (grado 6).

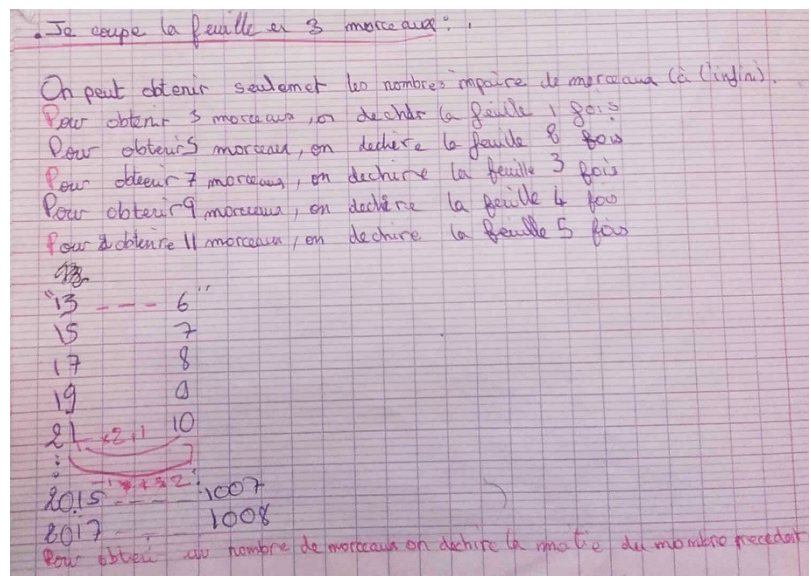

**Figura 5.** Evoluzione dell'esperienza su oggetti concreti e poi astratti.

Si nota in un primo momento una ricerca empirica, quando si taglia il foglio in tre pezzi ottenendo il risultato «Possiamo ottenere solamente numeri dispari di pezzi (all'infinito)» (seconda riga in **Figura 5**). In seguito, questa esperienza è riportata parola per parola fino a «Per ottenere 11 pezzi…» (settima riga). In seguito, inizia la formalizzazione matematica: non si parla più di pezzi di carta e di numeri di tagli, ma si fa una corrispondenza tra numeri:

- 13 — 6
- 15 — 7 ecc. fino a:
- 2017 — 1008.

---

[5] Per maggiori informazioni sulle analisi matematiche e didattiche di questo problema si rimanda a questo link: https://clarolineconnect.univ-lyon1.fr/icap_website/1324/40942

Infine, questa manipolazione degli oggetti matematici astratti sfocia in un'esplicitazione delle operazioni che legano le due colonne di numeri: $10 \times 2 + 1 = 21$ (esempio generico che si applica a tutti i numeri) e nell'altro senso: $-1 + \div 2$, che mostra ancora una difficoltà nel formalizzare la divisione euclidea, dando però il metodo che viene applicato per 2015 e 2017.

Questo esempio è rappresentativo del comportamento degli allievi confrontati a un problema di questo tipo e mostra bene le tre diverse manipolazioni di cui abbiamo parlato precedentemente: manipolazione dell'astratto che arriva progressivamente attraverso un graduale abbandono dell'esperienza concreta che serve a descriverla, come possiamo anche vedere in **Figura 6**; manipolazione dei simboli che permettono di esprimere le osservazioni emerse dall'esperienza, in particolare il modo di far corrispondere le colonne (**Figura 5**) e un tentativo di definizione per ricorrenza dei numeri della seconda colonna, utilizzando le differenze finite (**Figura 6**); manipolazione della logica, quando la validazione dei risultati emersi dalle manipolazioni sugli oggetti astratti si mette a confronto con i risultati dell'esperienza concreta (**Figura 4**).

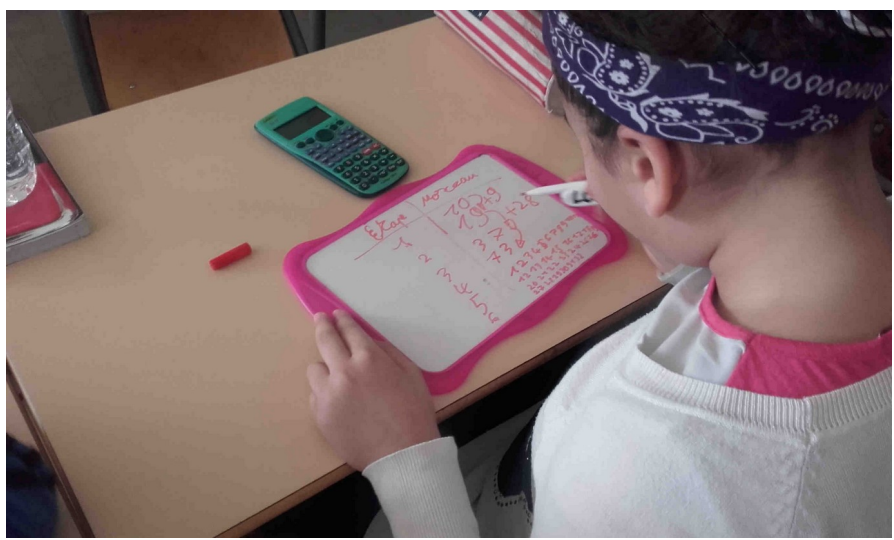

**Figura 6.** Passaggio a una manipolazione dell'astratto.

**4 Conclusioni**

La negoziazione di un contratto didattico nel quale gli allievi accettano che si possa apprendere la matematica, confrontandosi con dei problemi matematici,[6] non è per nulla scontata. Si tratta di un compito complesso che passa da una fase nella quale si stabilisce una cultura di ricerca all'interno della classe, attraverso dei quesiti, delle sfide, dei piccoli problemi, delle routine che possono allettare la curiosità e coinvolgere gli allievi in una nuova voglia di capire. Se la risoluzione di problemi diviene parte integrante del metodo di insegnamento e si fonda sulle manipolazioni (intese come descritte in questo articolo), essa diventa non solo un vettore di apprendimento di euristiche e di metodi di dimostrazione (competenze meta-matematiche) ma anche un modo di introdurre o di consolidare le conoscenze di oggetti matematici in legame con i piani di studio delle classi interessate. Gli esempi presentati in questo articolo mettono in evidenza gli apporti di questo approccio per gli apprendimenti matematici, mostrando l'importanza di differenti manipolazioni che sono necessarie per stabilire dei risultati matematici. Lungi dal contrapporre l'allenamento delle tecniche e la ricerca del senso, questi esempi mostrano bene come la costruzione dialettica della conoscenza matematica si appoggi invece sull'equilibrio tra manipolazione di oggetti astratti, manipolazione di simboli inseriti in un registro di rappresentazione dato, e manipolazione della logica (o delle logiche). Come esplicitato da Michèle Artigue a proposito del calcolo:

> «Ciò che costituisce la potenza della matematica, infine, non è solamente il fatto che essa si appropri di oggetti calcolabili e di sistemi di rappresentazione che sostengono efficacemente il calcolo, ma è anche il fatto che questo calcolo possa diventare algoritmo e automatizzarsi. Il calcolo è di conseguenza preso in un movimento altro, più potente, quello della sua meccanizzazione che, quando è riuscita, permette di svolgerlo senza pensare, riducendolo a una successione automatizzata di gesti. Questa meccanizzazione è necessaria per la progressione della conoscenza e vi è quindi, nella maggior parte dei calcoli, una sottile alchimia tra pensiero e routine».
>
> (Artigue, 2005, p. 4, traduzione dell'autore)

Le questioni che emergono sono dunque di ordine pedagogico e didattico e le si affronta attraverso sperimentazioni di un insegnamento basato sulla risoluzione di problemi, in cui la progressione annuale non sarà scandita in conoscenze da acquisire, ma piuttosto in problemi o in situazioni didattiche di problem solving. Tali situazioni lasciano la possibilità agli allievi di mobilitare le loro conoscenze e le loro competenze matematiche per fare emergere dei nuovi concetti o consolidare delle nozioni già presenti. L'ambizione dell'attuale lavoro dell'equipe DREAM è quella di esplorare su grande scala e in un contesto ordinario le condizioni e i vincoli per un insegnamento efficace della matematica, incentrato sul processo "manipolare–verbalizzare–astrarre" attraverso la risoluzione di problemi. I primi risultati mostrano da una parte la fattibilità di tale processo di insegnamento nelle classi "ordinarie" e dall'altra gli apporti in termini di apprendimento della matematica per gli allievi.

---

[6] Con l'espressione "problemi matematici" si intende ricordare che in questo articolo il focus è stato messo sul processo di matematizzazione verticale attivato dagli allievi, ossia sul lavoro da essi svolto all'interno del mondo della matematica per risolvere i problemi proposti.